\documentclass[12pt]{article}
\usepackage[utf8]{inputenc}
\usepackage{amsmath}
\usepackage{amsfonts,amssymb}
\usepackage{geometry}
\usepackage{mathtools}
\usepackage{amsthm,bm,todonotes}
\usepackage{authblk}
\usepackage{svg}
\usepackage[sorting=none]{biblatex}
\title{ \textbf{{Modelling to Generate Continuous Alternatives: Enabling Real-Time Feasible Portfolio Generation in Convex Planning Models}} }
\author[1,2,*]{Michael Lau}
\author[3]{Xin Wang}
\author[3]{Neha Patankar}
\author[1,2]{Jesse D. Jenkins}

\affil[1]{Andlinger Center for Energy and the Environment, Princeton University}
\affil[2]{Department of Mechanical and Aerospace Engineering, Princeton University}
\affil[3]{Department of Systems Science and Industrial Engineering, Binghamton University}
\affil[*]{Corresponding Author: Michael Lau, ml6802@princeton.edu}

\usepackage[T1]{fontenc}
\usepackage{fancyhdr}
\usepackage[english]{babel}
\usepackage{cancel}
\usepackage{tikz}
\usetikzlibrary{shapes.geometric}
\usepackage{float}
\restylefloat{table}
\usepackage{newclude}

\usepackage{xcolor}

\usepackage{graphicx}
\usepackage{subcaption}
\usepackage{caption}

\usepackage{float}
\usepackage{listings}
\usepackage{color}
\usepackage{textcomp}
\definecolor{dkgreen}{rgb}{0,0.6,0}
\definecolor{gray}{rgb}{0.5,0.5,0.5}
\definecolor{mauve}{rgb}{0.58,0,0.82}
\usepackage{fvextra}
\usepackage{csquotes}
\usepackage{setspace} 
\begin{document}
\maketitle
\newpage
\section{Introduction} \label{intro}
\paragraph{}{With the transformation of energy systems progressing worldwide, macro-energy systems (MES) models are increasingly valuable tools for providing decision support and information as to how decisions in the energy sector might impact outcomes of interest for a wide variety of stakeholders. Historically, most MES models have been tailored to address questions related to engineering feasibility, affordability and timing of transitions and are thus principally focused on the planning and operation of energy systems, namely cost-minimized capacity expansion and dispatch problems \cite{van_ouwerkerk_impacts_2022,jenkins_enhanced_2017}. }

\paragraph{ }{This set of least-cost optimization models work well for demonstrating the feasibility of a given scenario, providing a good estimate of the cost required to generate the desired outcome, and giving a way to understand complex dynamics present within the system. However, as argued by DeCarolis (\citeyear{decarolis_using_2011}) and Trutvenyte (\citeyear{trutnevyte_does_2016}), a least-cost optimization framework does not reflect the way decisions are made in the real world, nor can it allow for multiple stakeholders to balance their respective interests in planning outcomes. Feasible approaches to decarbonization can differ widely in their impact to both society at large and to individual stakeholders and interest groups \cite{sasse_distributional_2019,pedersen_using_2023}. Each stakeholder has their own objectives they would like to see reflected in the energy systems transition, the aims of which often appear to be in conflict. In real-world planning, while affordability is a key objective, it is one of several oft-competing priorities that must be balanced or satisfied, rather than optimized exclusively. }
\paragraph{ }{It is for precisely these reasons that modellers have looked to provide an understanding of the flexibility available in modelled systems, both in energy modelling and other sub-fields. Three primary tools have been used to expand the view of potential outcomes in optimization models. Parametric sensitivity analysis varies inputs to represent numerical uncertainty ranges for a subset of important parameters, which creates a set of scenarios with different assumption sets with the goal of testing how the model responds to parameter changes \cite{schyska_sensitivity_2021,yue_review_2018}. Multi-objective optimization (MOO) supplements the least-cost objective with a set of alternative objectives to generate a set of non-dominated solutions known as the Pareto frontier, within which no objective can be improved without sacrificing another, a condition called Pareto optimality or Pareto efficiency \cite{gunantara_review_2018}. Modelling to generate alternatives (MGA) imposes a budget above optimal cost, then reoptimizes repeatedly with varied objective functions to map the feasible space of near-optimal solutions that fall within the budget constraint \cite{decarolis_using_2011,trutnevyte_context-specific_2012,lau_measuring_2024}. }

\paragraph{ }{By presenting a highly diverse set of feasible near-optimal solutions, MGA allows stakeholders to decide their own relative priorities, shifting the epistemic burden of assuming stakeholder objectives away from researchers. This advantage comes with challenges. In addition to requiring substantial computational resources to compute enough MGA iterates to adequately approximate the feasible space of near-optimal solutions that fall within the budget constraint, MGA generates a lot of data - typically hundreds of fully fleshed-out MES model solutions for a given set of parametric assumptions \cite{lau_measuring_2024,pedersen_modeling_2021}. Thus, an important open research question remains: how can one best present MGA results in a coherent and accessible way that provides effective decision support in practical contexts? }
\paragraph{ }{MGA results are currently presented by either simplifying outputs or presenting more detailed results directly. There are two main methods of simplifying MGA outputs: a "common features" approach and sampling. The common features approach looks to identify "must-have", "must-avoid", and "real choice" characteristics between all solutions found \cite{voll_optimum_2015}. Common features can work well to tell coherent stories in papers, but this approach discards significant information in the process. Sampling is carried out by selecting a limited set of discrete solutions from within the set of near-optimal feasible solutions to present individually. The most commonly-employed algorithmic sampling method is to find the set of "most-diverse" solutions, that is the numerically-limited subset of MGA solutions that differ most from each other as first seen in Berntsen and Trutnevyte (\citeyear{berntsen_ensuring_2017}). Other papers frequently use heuristic sampling methods, where researchers select solutions which correspond to minimal or maximal values for certain outcomes of interests \cite{neumann_near-optimal_2021, chen_balancing_2022, sasse_distributional_2019}.  Especially with public-facing work, compiling all objectives of interest will be impossible, since each individual may have their own preferences. Even if comprehensive stakeholder elicitation is incorporated in the sampling process, previous research shows that stakeholders stated and revealed preferences are not identical and may change with new information \cite{pizarro-irizar_assessing_2020, trutnevyte_context-specific_2012}.}
\paragraph{ }{Previous work has also demonstrated that MGA outputs can be presented directly, either statically or interactively. Most MGA papers present MGA outputs statically, showing the range of outcomes for a given metric or variable found in the solution set \cite{price_modelling_2017,pedersen_modeling_2021}. However, such statistical measures may unintentionally mislead stakeholders or decision makers, as they can easily be interpreted as implying that the full range of outcomes as simultaneously achievable, while in reality, only individual specific combinations of decision variables are feasible and satisfy budget constraints (e.g. selecting a value from the range of outcomes in one dimension restricts the feasible set to a limited range of available values for other dimensions). Additionally, it is not clear that the frequency at which any given variable value appears within MGA iterates has any practical significance for decision support, as the appearance of an individual iterate within the near-optimal feasible set does not have any relationship to the real-world feasibility or likelihood of such an outcome.}
\paragraph{}{Pickering, Lombardi, \& Pfenninger (\citeyear{pickering_diversity_2022}) presented MGA solutions in an online dashboard which allows users to subselect individual MGA solutions based upon 9 different criteria and plots geographic breakdowns of primary energy consumption, solar and wind generation, transmission capacity expansion, and national imports and exports for each country in the EU. While interactive dashboards have been used in energy transition studies and stakeholder engagement before, as in Bessette et al. (\citeyear{bessette_decision_2014}), Mayer et al. (\citeyear{mayer_informed_2014}), and Volken, Xexakis, \& Trutnevyte (\citeyear{volken_perspectives_2018}), more research is needed as to how to make these interactive presentations effective communication, decision support, and education tools \cite{xexakis_are_2019}}. 

\paragraph{}{One commonality among the techniques described above is that they all present the discrete solutions identified by whatever MGA algorithm is used to explore the near-optimal feasible region. However, in the case of convex optimization problems, it is also possible to take advantage of the convexity of the space mapped by those solutions, as demonstrated by Pedersen et al. (\citeyear{pedersen_modeling_2021}). Pedersen et al. (\citeyear{pedersen_modeling_2021}) uniformly sample many potential solutions from the interior of the near-optimal feasible region by utilizing the properties of convex combination, with weights selected via random walk or Markov Chain processes and points selected by simplices created via Delauney triangulation \cite{pedersen_modeling_2021,pedersen_using_2023}.}

\paragraph{}{We, too, look to use the convexity of the MGA solution space to generate new feasible near-optimal solutions, taking advantage of the same mathematical principals as Pedersen et al. (\citeyear{pedersen_modeling_2021}) in the sampling stage of their Modelling All Alternatives (MAA) algorithm. However, their application specifies two things which ours does not. First, they specifically use the vertices of the simplices generated via Delauney triangulation within the original convex hull as the vertices with which to perform the convex combination \cite{pedersen_modeling_2021}. Additionally, they specify that the weights applied in the convex combination be generated at random \cite{pedersen_modeling_2021}. While this method is excellent for uniform sampling of the interior of the near-optimal feasible region, it runs into issues when we attempt to use it beyond this application or with high dimensional results. Finding simplex vertices via Delauney triangulation is computationally intractable for spaces of >10 dimensions with existing convex hull solvers, which rules out exploring spatially resolved models or models with more than 10 key decision variables of interest (e.g. power generation resources/technologies) through this process \cite{barber_quickhull_1996}. Furthermore, users have no control over exploration as the process can only return uniformly sampled points, not adapt to user priorities.}
\paragraph{}{In this paper, we present Modelling to Generate Continuous Alternatives (MGCA), an MGA post-processing methodology that can rapidly: a) create new MGA solutions with associated outcome metrics from anywhere within the interior or exterior of the feasible region previously identified via any MGA algorithm (i.e. the convex hull of previously identified MGA iterates), b) impose new objectives and constraints and rapidly calculate the impact of the optimal solution for this new objective or the remaining feasible set of near-optimal solutions after imposition of this new constraint, c) quickly approximate Pareto frontiers between a wide range of outcomes of interest, and d) retrieve capacity solutions for any point within the near-optimal feasible region to recreate more detailed outcomes using the full capacity expansion model. In Section \ref{Methods}, we describe the mathematical background and methodological contributions of MGCA. Section \ref{Results} will demonstrate several new capabilities the process creates. Section \ref{Conclusions} discusses potential use cases for MGCA and concludes.}

\section{Methods} \label{Methods}
\paragraph{}{In this Section, we describe the MGCA process and its mathematical foundations. First, we describe the Capacity Expansion Model problem in subsection \ref{CEMs}. We then highlight the four methodological contributions that underpin MGCA. We describe a method for MGA dimensionality reduction while retaining stakeholder-relevant information in subsection \ref{Process}. Then, in subsection \ref{FeasProb}, we introduce an exploratory problem that incorporates convex combinations in an optimization problem, allowing users to explore the MGA space in new ways.  In subsection \ref{MetricCalc}, we describe explicitly how the property of convex combinations can be used to calculate or retrieve metric values for different metric types. In subsection \ref{Retrieve}, we discuss a method for exporting MGCA solutions to a CEM to fully resolve portfolios with economic dispatch.}
\subsection{Capacity Expansion Models} \label{CEMs}
\paragraph{}{Energy system capacity expansion models (CEMs) are a class of linear or mixed-integer programs which optimize the planned capacity additions, retirements, and operation of an energy system. CEMs are thus generally formulated around capacity investment/retirement decision variables $\bm{x}$ and operational decision variables $\bm{y}$. They also incorporate a variety of operational and planning constraints, such as transmission constraints, unit commitment for thermal generators, and generation decisions in each time-step $t \in T$. CEMs are typically formulated as follows:
\begin{equation}
    \label{eq:LP}
    \begin{alignedat}{3}
        &\text{minimize}&\; \; &\bm{c}^T\bm{x}+\sum_{t \in T} \bm{d}_t^T \bm{y}_t\\
        &\text{subject to}& & \bm{A}_t\bm{x} + \bm{B}_t\bm{y}_t \leq \bm{b}_t, \quad \forall t \in T\\
        &&& \bm{y}_t \in Y_t,\quad \forall t \in T\\
        &&& \bm{x} \in X
    \end{alignedat}
\end{equation}
Here, $\bm{A}_t\in \mathbb{R}^{r \times n}$ and $\bm{B}_t\in \mathbb{R}^{r\times m}$ are constraint coefficient matricies, $\bm{b}_t \in \mathbb{R}^r$ the vector of right hand values, and $\bm{c}\in \mathbb{R}^n$ and  $\bm{d}_t\in \mathbb{R}^m$ represent fixed and operational costs. $X \subseteq \mathbb{R}^n$ and $Y_t\subseteq \mathbb{R}^m$ are polyhedral sets representing linear constraints on planning and operational decisions. In the MGA context, this formulation changes to the following: 
\begin{equation}
    \label{eq:MMGA}
    \begin{alignedat}{3}
        &\text{minimize}&\; \; &\bm{w}^T\bm{x}\\
        &\text{subject to}& & \bm{c}^T\bm{x}+\sum_{t \in T} \bm{d}_t^T \bm{y}_t \leq (1+\epsilon)\Big(\bm{c}^T\bm{x}^* + \sum_{t \in T}{\bm{d}_t^T\bm{y}_t^*}\Big) \\
        &&& \bm{A}_t\bm{x} + \bm{B}_t\bm{y}_t \leq \bm{b}_t, \quad \forall t \in T\\
        &&& \bm{y}_t \in Y_t,\quad \forall t \in T\\
        &&& \bm{x} \in X
    \end{alignedat}
\end{equation}
Where $\bm{w} \in \mathbb{R}^n$ is the MGA coefficient vector and $\epsilon > 0$ is the user-defined budget slack. } 
\subsection{Process Description and Dimensionality Reduction}  \label{Process}
\begin{figure}[H]
    \centering
    \includegraphics[scale=0.6]{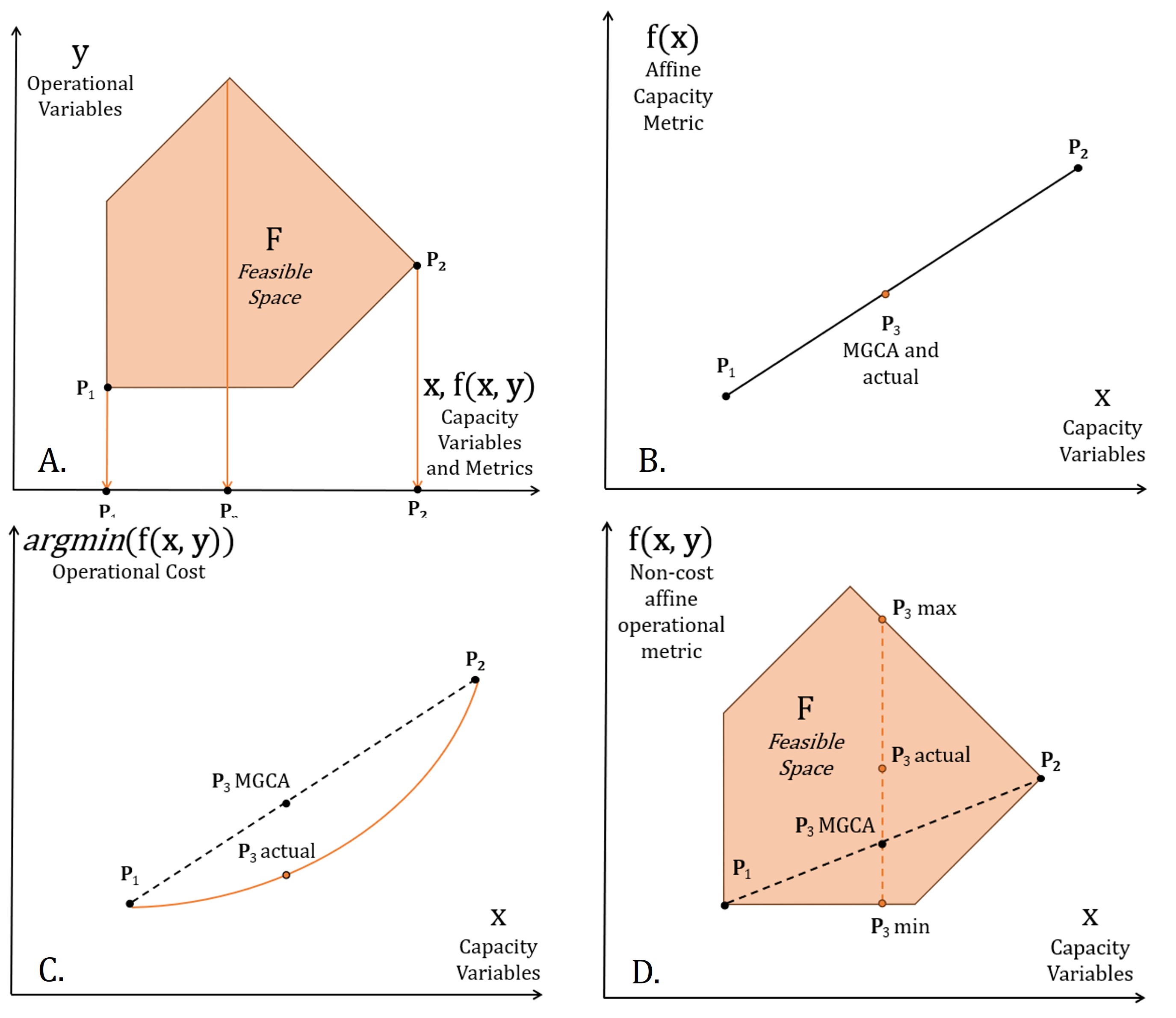}
    \caption{MGCA dimensional reduction and metric interpolation visualization. A. Dimensional reduction from capacity and operational variables and metrics to capacity variables and metric values. B. Convex combination of affine capacity metrics, which can be precisely calculated. C. Convex combination of operational cost, the actual value of which will always be less than or equal to the interpolated value, as it is an \textit{argmin}, a convex function. D. Interpolation of other operational metrics, which is guaranteed to produce a feasible value, but not the actual value.}
    \label{fig:MGCALogic}
\end{figure}

\paragraph{}{Modelling to Generate Continuous Alternatives utilizes the properties of convex combination to rapidly generate new, feasible, within-budget solutions to the LP in response to stakeholder priorities. We do this through a series of steps. First, one or more MGA algorithms are run on a given problem to generate a set of solutions which explore as much of the near-optimal feasible space as possible (see \textcite{lau_measuring_2024} for more on methods to efficiently and extensively explore the near-optimal feasible region of an MGA problem).} 
\paragraph{}{We then reduce the dimensionality of the original MGA solutions by projecting the set of MGA solutions onto the subspace formed by a set of capacity decisions and selected metrics (e.g. outcomes of interest). Specifically, we create a new subspace $Z$ which contains a selected set of capacity decisions from all  solutions, or iterates, identified by MGA, $I$, and a set of metrics of interest which are linear functions of decision variables, $M$. Formally, $(I \cup M) \subseteq Z$. Thus, the elements of $Z$ are capacity decisions $\bm{x}$ and metrics of interest as functions of operational and capacity decisions $f^T \bm{y}$ and $g^T \bm{x}$ respectively where $f\in \mathbb{R}^m$ and $g \in \mathbb{R}^n$ are coefficient vectors defining each metric.}  
\paragraph{}{Dimensionality reduction reduces solve times for all operations on the MGA points by reducing the length of each point vector, but must retain most stakeholder-relevant information. In a CEM, there are $O(N)$ capacity decisions, where $N$ is the number of existing or candidate generation, storage, and transmission capacity options, and $O(N\times T)$ operational decisions where $T$ is the number of operational time-steps, usually $O(10^3)$ to $O(10^5)$. The dimensionality reduction proposed here ensures that MGCA solves problems of size $O(M\subseteq N)$ rather than of $O(N + (N \times T))$ as would occur if applied to a full capacity expansion problem, a 3-5 order of magnitude reduction. In a policy-facing setting, the goal of MGA is typically to communicate to users a set of alternatives which map the flexibility in investment and capacity decisions possible within a given cost slack, as well as the implications of those decisions in the form of metrics for various outcomes of interest, such as implications for land use, air pollution, employment, tax revenue, energy import dependence, etc. \cite{lombardi_policy_2020,pedersen_using_2023,sasse_distributional_2019,sasse_regional_2020}. As such, we retain capacity decisions and affine metric values of interest, including but not limited to cost, emissions, and land use both for the system and specific zones, while discarding operational decisions as shown in panel A of Figure  \ref{fig:MGCALogic}. These lower-dimension points can then be used as vertices for the convex hull of the near-optimal feasible region of the planning problem, the interior of which can be explored via convex combination. }
\subsection{Convex Hull Exploration Problem} \label{FeasProb}
\paragraph{}{Linear Programming models are created through the intersection of a set of linear halfspaces, defined by the inequality constraints of the model, which form a convex polyhedron called the feasible space \cite{bertsimas_introduction_1997}. The feasible space is by necessity convex and affine as all constraints are linear \cite{bertsimas_introduction_1997}. }
\paragraph{}{We can leverage the properties of this geometry through the introduction of convex combinations. The convex combination of points is defined as the weighted sum of a set of points wherein all weights are positive and the sum of the weights is equal to one \cite{bertsimas_introduction_1997}. This can be written formally as:
\begin{gather}
    \bm{z}_1,\bm{z}_2,\bm{z}_3,\dots,\bm{z}_j \in Z \\ 
    \bm{\lambda}_i \geq 0, \sum_{i=1}^{j}\bm{\lambda}_i = 1 \\
    \sum_{i=1}^{j}{\bm{\lambda}_i\bm{z}_i} = \bm{z}_k \in Z \label{eq 1}
\end{gather}
where $Z$ is a polyhedral set, points $\bm{z}_i\in \mathbb{R}^h$, and $\bm{\lambda}\in \mathbb{R}^h$ the set of weights of the convex combination. The convex hull of a set of vectors, meanwhile, is defined as the "set of all convex combinations of these vectors."\cite{bertsimas_introduction_1997} Thus, by definition, we can recreate any point within the convex hull using only the convex combinations of vertices of the convex hull.}
\paragraph{}{To enable solutions within convex hull of the near-optimal feasible space to be explored using the breadth of existing optimization tools, we introduce an exploratory optimization problem, shown in equation \ref{feaspr}, for creating weights in response to user priorities, constraints, objectives, and perceptions. 
\begin{equation}
    \begin{aligned}
    \label{feaspr}
        min \:&f(\bm{z}_k) \\
        s.t. \: &\bm{z}_k = \bm{\lambda}^T Z\\
            &\sum{\bm{\lambda}} = 1\\
            &\bm{Az}_k \leq \bm{b}\\
            &\bm{\lambda} \geq 0\\
            &\bm{z}_k \geq 0
    \end{aligned}
\end{equation}
Here $\bm{z} \in \mathbb{R} ^ n$ denotes the vector of $n$ capacity decisions and metrics generated by the equation, $Z \in \mathbb{R}^{n \times m}$ denotes the set of $m$ MGA vertices (the dimensionality-reduced $Z$ from Section \ref{Process}, $\bm{\lambda} \in \mathbb{R}^m$ indicates the vector of weights, which are decision variables in this model.  The program functions as follows. The objective may be specified by the user. It could sample the space by repeatedly solving for different vectors, created via MGA vector selection (see \textcite{lau_measuring_2024}), or solve for other specified objectives or metrics the user is interested in. The constraints in equation \ref{feaspr} constrain generated solutions $\bm{z}$ to be within the set of convex combinations defined by $Z$ and $\bm{\lambda}$, which guarantees any solutions to \ref{feaspr} will fall within the within-budget feasible set of the original MGA problem \ref{eq:MMGA}. Coefficients $\bm{A} \in \mathbb{R}^{n 
 \times m}$ and $\bm{b} \in \mathbb{R}^{m} $ define generic constraints on the capacity decisions and any metrics that are affine products of capacity decisions in $Z$ which may be added by the user. The exploration problem has $2m$ variables and $3n+d+1$  constraints, where $d$ is the number of user-added constraints in $\bm{Az}_k \leq \bm{b}$. This is typically $O(10^{3}-10^4)$ and therefore can be evaluated live. Note that due to the small size of the problem, it is entirely possible to use non-linear convex objectives. Furthermore, since the exploration problem is an optimization problem, it will, by definition, always find an exterior solution. Unless new constraints are added, this means that it will always find an original MGA vertex. One can thus explore the interior of the space by either specifying weights manually or by imposing additional constraints not found in the original MGA problem \ref{eq:MMGA} then optimizing.}
 \paragraph{}{The power of the convex hull exploration problem lies in that it allows users to explicitly identify optimal solutions within the captured near-optimal feasible region accounting for additional user-specified constraints and/or alternative user-specified objectives. As an example, suppose a utility is engaged in an integrated resource planning process to determine capacity expansion plans and is discussing potential capacity plans with local groups using an MGCA based analysis. Perhaps a local environmental group wants to minimize the amount of prime farmland which will be converted to renewables, without restricting it entirely. They could see how much prime farmland \textit{must} be converted to renewables without increasing total costs by more than a specified amount (the MGA budget slack constraint) by solving the exploration problem with a new objective minimizing the sum of renewable sites on prime farmland, thus giving them a beginning bargaining position and enabling them to demonstrate the feasibility of their request. A climate group may want the same system to limit the region's build-out of new natural gas generators. They could run the exploration problem on a least-cost system with a set of increasingly stringent constraints on new gas capacity additions to see how that impacts least-cost system build-out. The two groups could then negotiate or work together by imposing the constraint on the prime farmland-minimizing problem to see how their desires impact one another. One could imagine multi-party negotiations occurring this way, with groups finding solutions that suit them and combining them with other groups preferred solutions as a means of compromise. As solutions to problem \ref{feaspr} are calculated in seconds or less, this method could also permit the design of interactive decision support tools that allow various stakeholders and decision markets to quantitatively explore the full range of near-optimal feasible solutions to a capacity expansion planning problem and build intuition as to the shape of important trade-offs and considerations.}
\subsection{Metric Calculations}\label{MetricCalc}
\paragraph{}{We can further utilize the properties of convex combinations and the scalar product to evaluate affine metrics for each newly constructed feasible point. By definition, the scalar product is scalar-multiplicative \cite{lipschutz_vector_2009}. Thus, we can take the convex combination of an affine metric of the vertices of the convex hull to calculate the result for any point within the space. Formally, we can write this as: }
\begin{gather}
c^\intercal(\sum_{i=1}^{n}{\lambda_ix_i}) = \sum_{i=1}^{n}{\lambda_i(c_ix_i)}
\end{gather}
\paragraph{}{When this property is applied to the capacity and metric projection of the original set of MGA solutions created by the first step of MGCA, we find that affine metrics of capacity variables can be evaluated to exact values in post-processing, regardless of if they were initially evaluated in the model, as shown in panel B of Figure \ref{fig:MGCALogic}, because all relevant capacity data is retained within the dimensionality reduction process. }
\paragraph{}{We can also use this property to estimate values of metrics that are an affine product of operational decision variables and computed in the initial MGA run, like operational cost and emissions, through convex combination. As MGCA simply interpolates between metric values rather than re-optimizing the operations of the new system, created metric values will not precisely match what would occur under economic dispatch. Since operational cost is minimized by economic dispatch, MGCA estimates of system-level variable cost will always be greater than or equal to  the system cost as evaluated by a least-cost operational model, given the same assumptions as shown in  Figure \ref{fig:MGCALogic}, panel C. However, this is not true of zonal or regional costs, which may increase in specific zones relative to the MGCA estimate in order to bring down system operational costs. Since other metrics are not present in the dispatch objective problem, no guarantee of lower overall system metric values is possible. Convex combination of operational metric values will not match economic dispatch, but will be feasible values within the convex hull of those found through the full model with economic dispatch, by definition as shown in panel D.  However, it is possible to constrain other metrics through policy constraints, like an emissions cap, which will constrain the feasible space, and thus the MGCA interpolates, to guarantee a certain maximum level of emissions.}
\subsection{Exporting MGCA Solutions for Full Computation}\label{Retrieve}
\paragraph{}{Once users use the exploration problem to find solutions they are interested in, they may wish to explore a fully fleshed-out solution with similar attributes and technologies emphasized. Computing a full solution with least-cost dispatch, unit commitment, and other key details provides assurance that the discovered solution is in-fact feasible, and provides accurate evaluations of operational metrics consistent with least-cost economic dispatch. We thus provide a method to reproduce discovered solutions in a full CEM.}
\paragraph{}{It is possible to replicate the capacity decisions in an MGCA interpolate within a full CEM by fixing all capacity decisions in the CEM to those in the interpolate and solving the CEM with a least-cost objective. By fixing the capacity variables, we reduce the complexity of the model to an operational model which shortens runtime while preserving least-cost operations. It is worth noting that since the capacities of the CEM solution are fixed with equality constraints, their values will not be optimized, resulting in a total cost below the MGA budget imposed on the initial problem as the solution is in the interior region of the capacity space. This method can recreate the capacities of any MGCA solution with least-cost operational decisions, either in the interior or on the exterior of the space.}
\section{Results }\label{Results}
\paragraph{}{In this Section we demonstrate three potential uses for the mathematical constructs discussed in Section \ref{Methods}.  We first discuss testing methods in subsection \ref{Testing}, covering the demonstration problem on which we show MGCA capabilities. In subsection \ref{explore} we demonstrate the use of the feasibility problem to explore the extent of the discovered MGA space for metrics of interest. Then in subsection \ref{Omit}, we demonstrate the addition of user-specified constraints to the space. In subsection \ref{pareto}, we demonstrate how MGCA can be used to approximate Pareto frontiers between outcomes of interest. Finally, in subsection \ref{ResultsRetrieve}, we test and discuss the accuracy of MGCA estimates of operations-based metrics. }
\subsection{Testing} \label{Testing}
\subsubsection{Demonstration Problem}\label{DemoMethod}
\begin{figure}[H]
    \centering
    \includegraphics[scale=0.5]{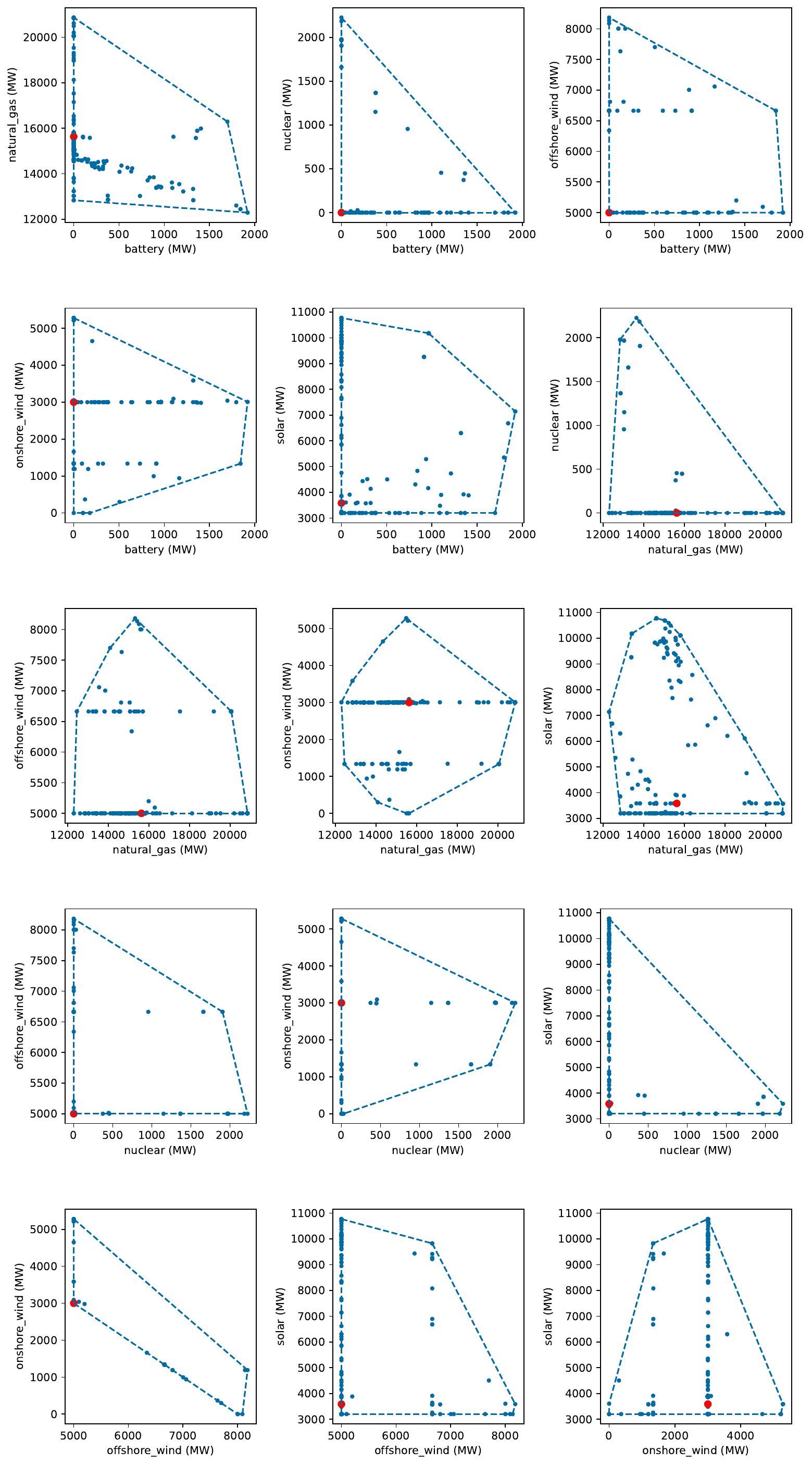}
    \caption{Aggregate capacities for all MGA solutions found for 3-Zone GenX ISONE Case. Solutions indicated by dots, convex hull indicated by dashed line. Least-cost solution in red.}
    \label{fig: BaseMGA}
\end{figure}

\paragraph{}{To demonstrate the capabilities provided by MGCA in a realistic setting, we computed an example capacity expansion model result. We used the open-source GenX capacity expansion model\cite{jenkins_genx_2023}, utilizing the three-zone ISONE example case in the GenX repository with 8760-hourly resolution and linearized unit commitment for a single year.  We ran MGA using the Random Vector method applied to capacity variables with a 10\% budget slack and computed 200 total MGA iterates, not including the least cost solution \cite{berntsen_ensuring_2017,lau_measuring_2024}. Problems were solved using Gurobi \cite{gurobi_optimization_llc_gurobi_2024}. The capacity decisions and operational metrics of the resulting MGA solutions have been aggregated into 12 single technology metrics representing all capacity for each technology in all zones and 2 metrics for the whole system, cost and emissions. We only visualize pairwise convex hulls of 6 technologies to reduce the number of subplots required. Please note, we aggregate results for simplicity of static visualization, not due to any computational restriction of MGCA (that is, one could retain individual capacity decisions for each resource in each zone or aggregate). The base MGA solutions found by this exploration are displayed in Figure \ref{fig: BaseMGA}.}
\paragraph{}{All MGCA post-processing computations were carried out on a laptop with an Intel i5 processor using the free solver GLPK\cite{oki_basics_2012}.}
\subsection{Feasible Space Exploration} \label{explore}
\begin{figure}[H]
    \centering
    \includegraphics[scale=0.5]{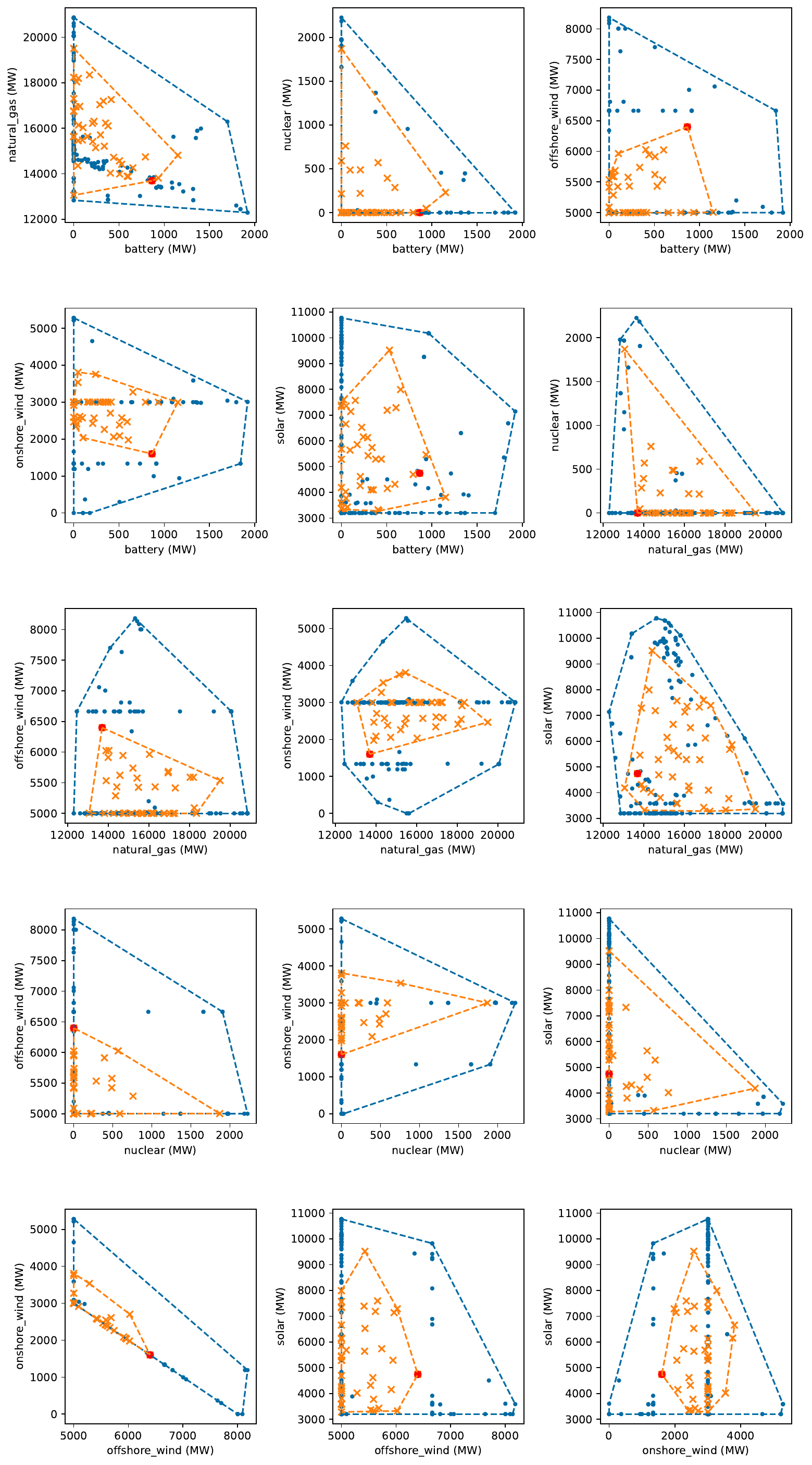}
    \caption{Capacities of 50 randomly generated interior MGCA interpolates (orange) and original MGA points (blue) }
    \label{fig:MGCAInterps}
\end{figure}

\begin{figure}[H]
    \centering
    \includegraphics[scale=0.5]{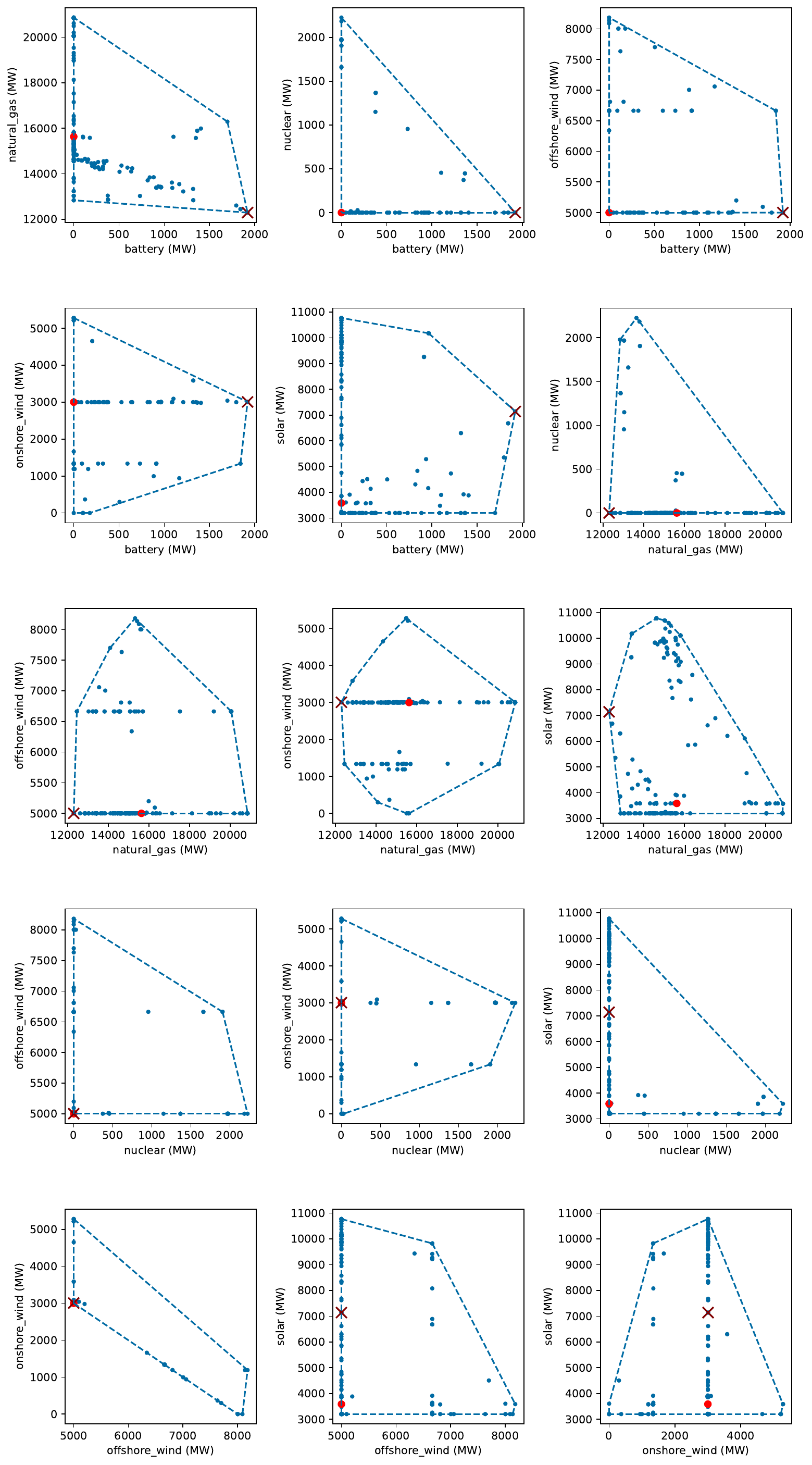}
    \caption{Capacities of solution found through exploration problem with objective minimizing natural gas (maroon X) and original MGA points (blue). Least-cost solution in red.  }
    \label{fig:MGCAFeasProb}
\end{figure}

\paragraph{} {Applying MGCA to model results sets, like those generated by MGCA and connecting them to an interactive dashboard allows users to generate new solutions not computed by the initial computation, including those anywhere in the interior of the near-optimal feasible region (without resorting to computationally expensive uniform sampling of the interior region as in MAA). We have considered several ways users can generate new solutions. In a case where users find solutions they like and want to find new, similar solutions, they can specify the sets of points they wish to interpolate between and their desired weights, which will rapidly generate the number of new solutions specified. As discussed in Section \ref{MetricCalc}, since we can evaluate affine functions of capacity variables to exact values and affine functions of generation variables to feasible values and generate these as part of the new point created by the user, providing them a way to rapidly evaluate how different capacity choices behave in the system being studied. This method is well suited for users who desire a large number of new points with a good awareness of which points they like and how they wish to weight them. }
\paragraph{}{We show an example of this capability by creating 50 MGCA interpolates, each with random weightings of 4 exterior points. The interior points were computed in 0.03 seconds in total, or .6 ms each, while the full MGA problem typically takes around 600 seconds to compute each iterate. Capacities for all points are shown in Figure \ref{fig:MGCAInterps}.}
\paragraph{} {Alternately, users can input search vectors within the space, such as finding the portfolio with the maximal solar in a given zone, the least wind in another zone, or the lowest total emissions. Using the exploratory problem presented in Section \ref{FeasProb}, we can rapidly find the initial MGA solution that optimizes the desired objective. This is analagous to traditional Multi-Objective Optimization, but here the optimization is performed in 'real-time' by solving a much smaller, pre-mapped decision space rather than on the entirety of the optimization model. If these vectors of interest are known prior to the initial computation, it is best practice to explicitly optimize for them as part of the MGA analysis as a bracketing run to ensure that the set of MGA iterates includes the maximum and/or minimum feasible values for this outcome of interest (see Lau, Patankar \& Jenkins, 2024). Without this step, the feasibility problem will discover the closest extreme point within the MGCA subspace to the optimal value in the larger model decision space. We demonstrate this capability in Figure \ref{fig:MGCAFeasProb} by evaluating the exploratory problem with a new objective: minimizing natural gas capacity in the system. Solving the exploratory problem returns the best MGA iterate for the objective, and in this case took 0.02 seconds, compared to around 600 to solve for a similar point with a full CEM optimization.}
\subsection{Evaluating Results Under Different Budgets} \label{Budgets}
\begin{figure}[H]
    \centering
    \includegraphics[scale=0.5]{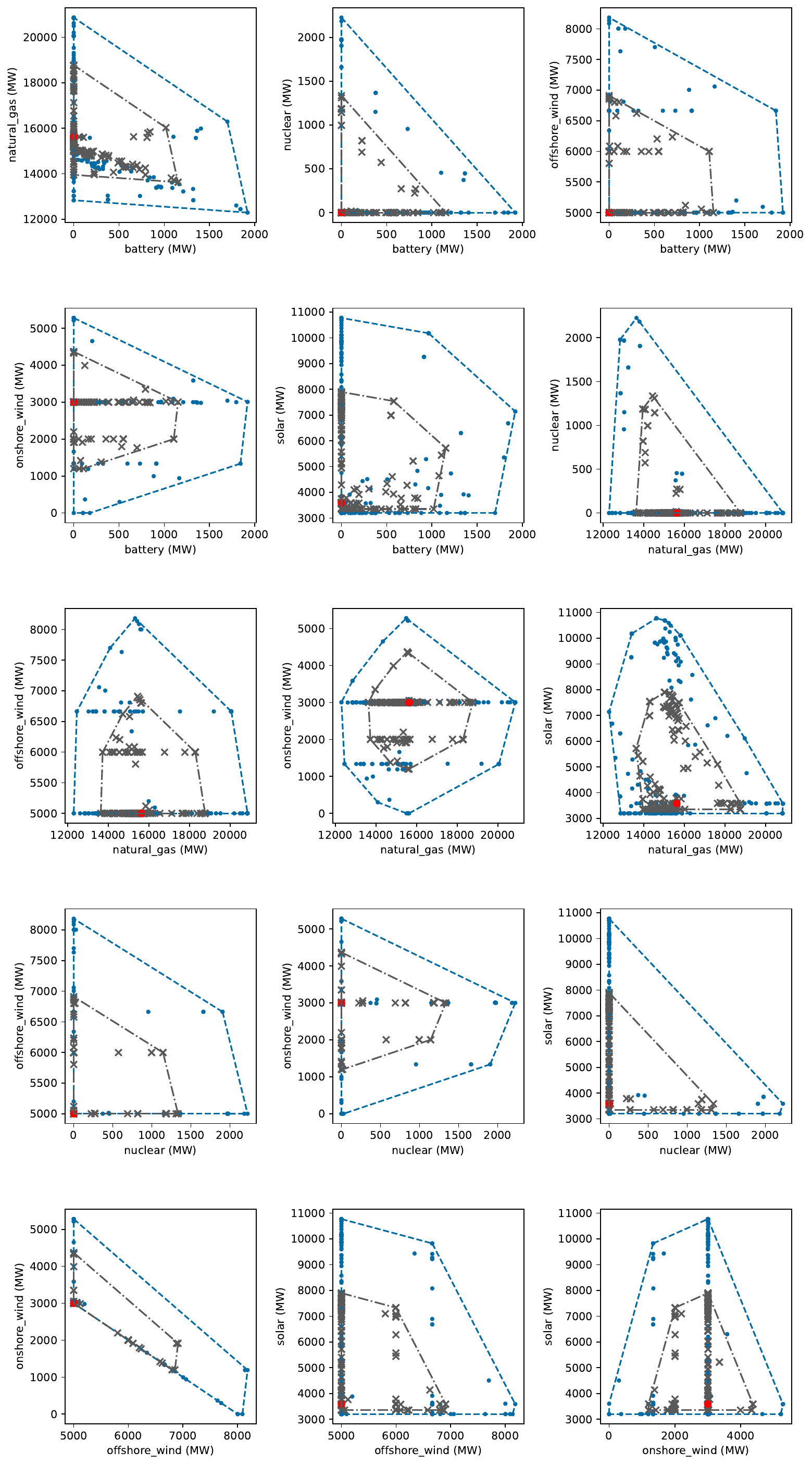}
    \caption{MGA capacity space with 6\% budget (gray) superimposed over 10\% budget (blue). Individual solutions are shown by x's (6\%) and dots (10\%) respectively. Least-cost solution in red.}
    \label{fig:BudgetCap}
\end{figure}

\begin{figure}[H]
    \centering
    \includegraphics[width=\linewidth]{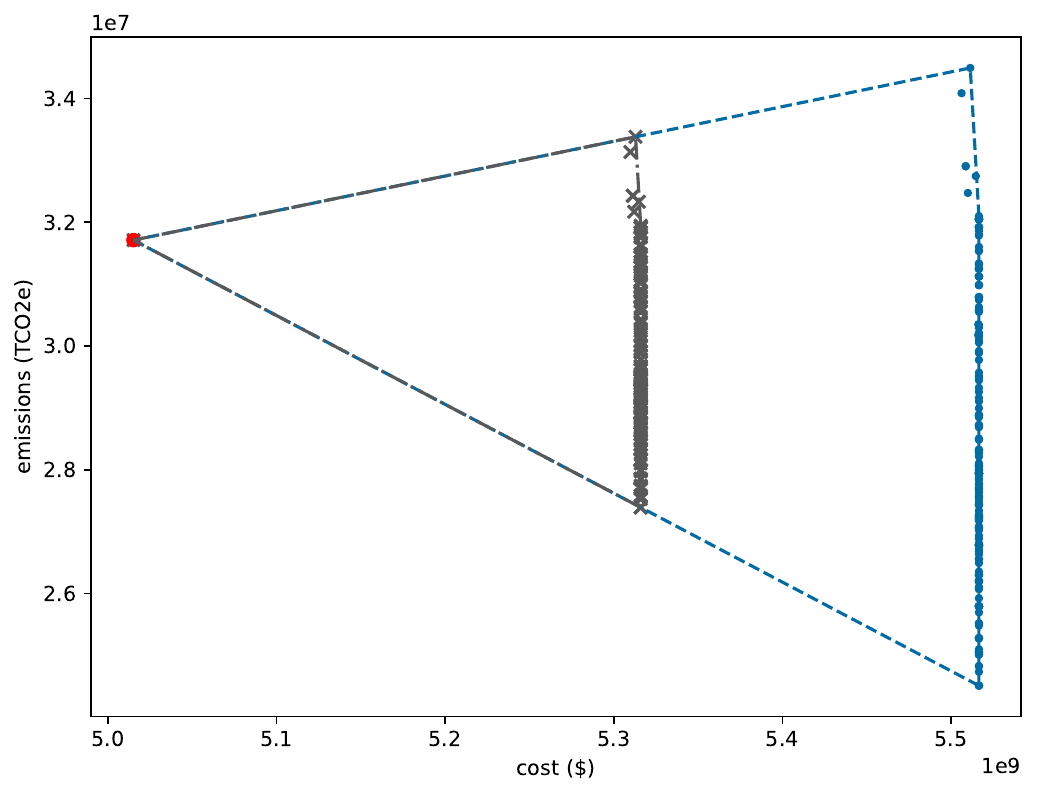}
    \caption{MGA metric space with 6\% budget interpolates (gray) superimposed over 10\% budget MGA points (blue).  Individual solutions are shown by x's (6\%) and dots (10\%) respectively. Least-cost solution in red.}
    \label{fig:budgetmetric}
\end{figure}

\paragraph{}{One potential application of MGCA is to generate MGA results for varying budgets. As demonstrated by \textcite{neumann_near-optimal_2021}, the degree of flexibility available at different budget levels can be salient to decision-makers. Under existing methods, a new set of MGA iterates are computed for each new budget level, incurring significant computational cost (as in \cite{neumann_near-optimal_2021}). With MGCA, it is possible to directly interpolate budget-constrained solutions extremely rapidly, and to represent the convex hull of the feasible region within this tighter constraint, effectively performing an entirely new MGA problem without having to rerun the CEM again. The steps to do so are listed below.}
\begin{enumerate}
    \item Run initial MGA at highest desired budget level
    \item Create MGCA projection onto capacity and metric subspace
    \item Iteratively combine all MGA solutions individually with least-cost solution, with weights specified by the proportion of the budget desired at each cost level. For example, with an initial budget of 10\%, an MGA solution at a budget of 6\% could be interpolated by weighting said point by 0.6 and the least-cost solution by 0.4 then summing them.
\end{enumerate}
\paragraph{}{The resulting points will populate the desired budget constraint. We demonstrate this capability in Figure \ref{fig:BudgetCap} }
\paragraph{}{Additionally, in Figure \ref{fig:budgetmetric} we plot the cost/emissions trade-off figure to demonstrate that the generated solutions do, in fact, relate to a 6\% budget and span the full range of MGA iterates found. The budget interpolation of all solutions computed in 0.35 seconds.}
\subsection{Incorporating Additional Constraints} \label{Omit}
\begin{figure}[H]
    \centering
    \includegraphics[scale=0.5]{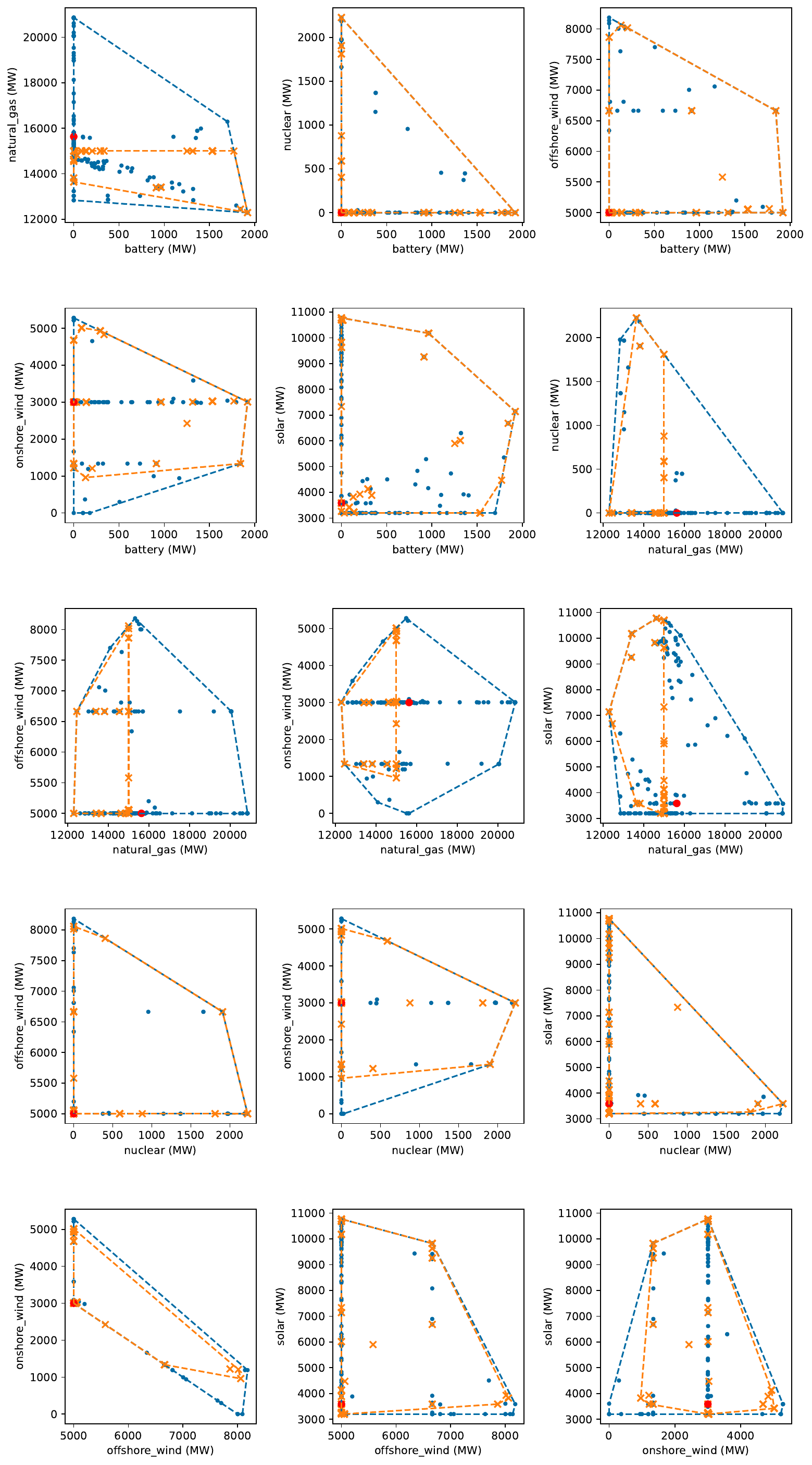}
    \caption{MGA capacity space with constraint keeping natural gas capacity $\leq$ 15000 MW added (orange) superimposed over original MGA space (blue). Individual solutions are shown by x's (constrained space) and dots (unconstrained) respectively. Least-cost solution in red.}
    \label{fig:Constraint}
\end{figure}

\paragraph{} {Constraints on the MGCA subspace can also be evaluated live. Since we have formulated a new feasibility problem based on the solutions found, we can add new constraints on the variables present in the feasibility problem and compute new extreme points created by the further constraint. Once constraints have been added, it is possible to investigate the impact of newly added constraints on the feasible space through repeated optimization using MGA vectors, similarly to ordinary MGA, just without a new budget constraint, as the polyhedron formed by the set of MGA solutions $Z$ already is constrained to the original budget. Carrying out a sufficiently large number of these exploratory optimizations will eventually approximate the full shape of the newly constrained space. Since there are only typically a few hundred variables in the exploratory problem, these optimizations can be carried out extremely quickly on personal computers.}
\paragraph{}{To demonstrate this capability, we added a new constraint that constrains natural gas capacity to 15000 MW to the feasibility problem and carried out 200 local MGA iterations using the Variable Min/Max method to explore the newly added constraint \cite{trutnevyte_context-specific_2012}. It is worth noting that due to the dimensional reduction carried out at the beginning of MGCA, this process of 200 MGA solutions including the least-cost, took 2.17 seconds or about 0.01 seconds per solution. The results are shown in Figure \ref{fig:Constraint}. This approach can be used in decision support contexts by permitting stakeholders or decision makers to insert additional user-defined constraints or preferences on the original planning problem, and then quickly visualize how this additional constraint (or constraints) shrink the remaining set of feasible solutions within the original budget constraint. This process could be repeatedly successively as a user adds new constraints or loosens/tightens constraints to build intuition about how these solutions limit the range of feasible planning options. }
\subsection{Creating Pseudo-Pareto Frontiers} \label{pareto}
\begin{figure}[H]
    \centering
    \includegraphics[scale=0.5]{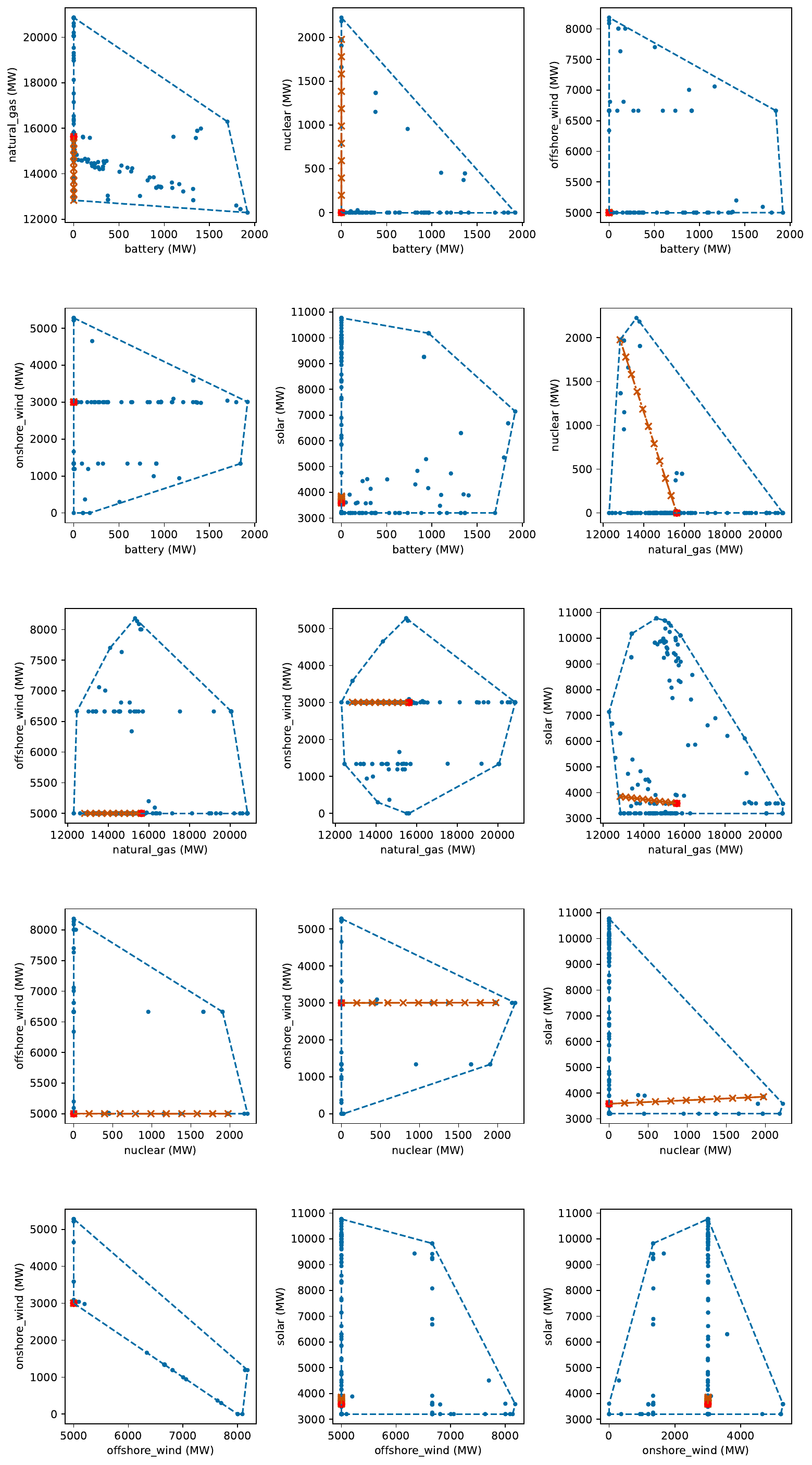}
    \caption{MGA capacity space with cost-emission Pareto frontier (set of non-dominated solutions between cost and emissions) (burnt orange) superimposed over original MGA space (blue). Individual solutions are shown by burnt orange x's (Pareto frontier) and blue dots (MGA space) respectively. Least-cost solution is highlighted in red.}
    \label{fig:ParetoCap}
\end{figure}

\begin{figure}[H]
    \centering
    \includegraphics[width=\linewidth]{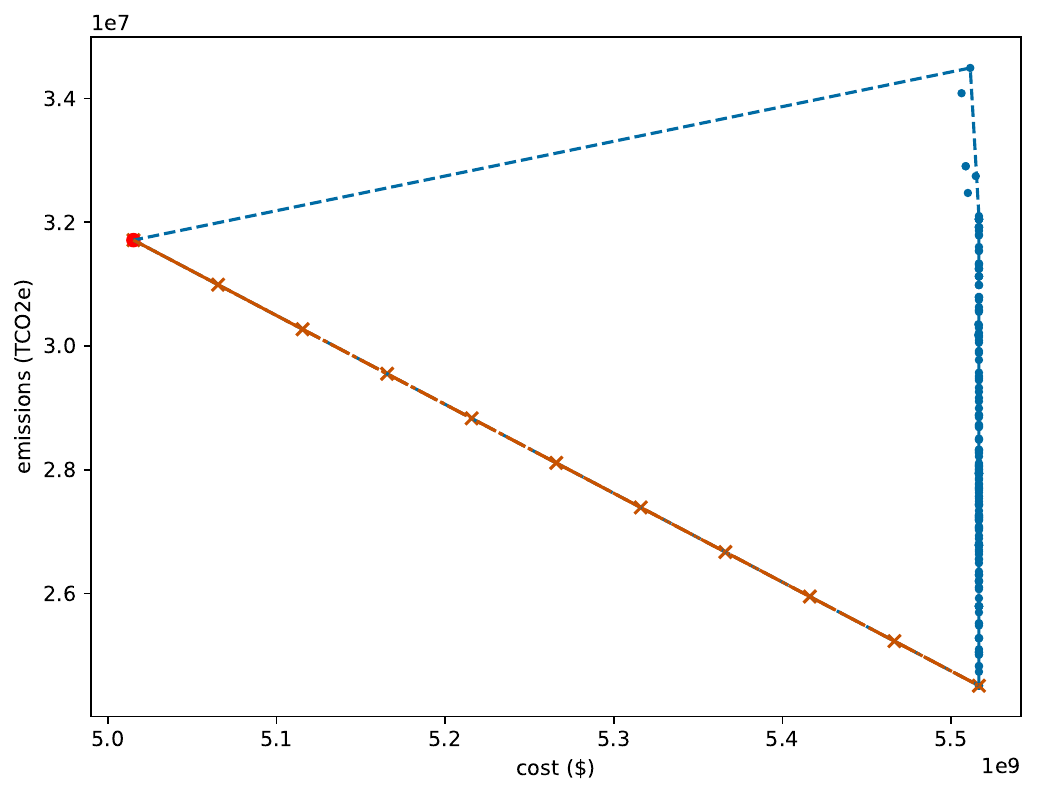}
    \caption{MGA metric space with Pareto frontier interpolates (orange) superimposed over original MGA points (blue).  Individual solutions are shown by burnt orange x's (Pareto frontier) and dots (MGA) respectively. Least-cost solution is highlighted in red.}
    \label{fig:ParetoMetric}
\end{figure}

\paragraph{} {Using MGCA, we can approximate Pareto frontiers between selected MGA solutions by running Multi-Objective Optimization (MOO) using a standard epsilon constraint method \cite{gunantara_review_2018,mesquita-cunha_new_2023} within the capacity/metric subspace using the feasibility problem in Section \ref{FeasProb}.  To create a Pareto frontier, the user selects two or more objectives, then a function would allow the feasibility problem to repeatedly find the extreme points representing the Pareto frontier of those objectives. The mathematical formulation of an MOO-type MGCA feasibility problem is shown in equation \ref{MOO}. This Pareto frontier is may not be identical to a similar Pareto frontier calculated on the original space, as it is constrained by the convex hull of computed MGA iterates (including the budget constraint) and is likely to not include the true optimal solution for most combinations of objective and epsilon constraint. However, this approach will find the optimal points within the capacity/metric subspace for the set of constraints and objectives selected and the specified MGA budget slack.}
\begin{equation}
\begin{aligned}
\label{MOO}
    min \:&g(\bm{z}) \\
    s.t. \:  &f(\bm{z}) \leq f(\bm{z}^{\ast})*(1+\epsilon)\\
        &\bm{z} = \bm{\lambda}^T Z\\
        &\sum{\bm{\lambda}} = 1\\
        &\bm{\lambda} \geq 0\\
        &\bm{z} \geq 0
\end{aligned}
\end{equation}
\paragraph{}{In Figure \ref{fig:ParetoCap}, we demonstrate this capability  by using MGCA to create a Pareto front between total system cost and total system emissions, two metrics of interest computed in our test case. The solution of these 11 problems took 0.09 seconds, or approximately 0.008 seconds per problem.}
\paragraph{}{By plotting the cost-emissions trade-off space, as shown in Figure \ref{fig:ParetoMetric}, we can see that the MGCA-generated solutions do in fact trace the trade-off between cost and emissions within our MGA space.}

\subsection{MGCA Operational Metric Accuracy} \label{ResultsRetrieve}
\paragraph{}{As discussed in Section \ref{MetricCalc}, MGCA can rapidly generate exact values of affine capacity-related metrics, but only produces \textit{feasible} values of generation-based metrics. It is important to understand how different feasible generation metric values are from the values given by a full CEM for a given problem in interpreting MGCA findings. }
\paragraph{}{One way to consider the accuracy of MGCA metric estimates is percent difference: how far away are MGCA metric estimates from the CEM calculated values with hourly economic dispatch and linearized unit commitment? To evaluate this, we took the 50 interpolates generated in Section \ref{explore} then exported all of them to the three-zone demonstration problem in GenX. We record cost and emissions values for each MGCA interpolate and their twin in the demonstration problem. In Figure \ref{fig:PercentDiff}, we see the percent difference between MGCA's estimated metric values and the actual values calculated through the full CEM. 
\begin{figure}[H]
    \centering
    \includegraphics[width=\linewidth]{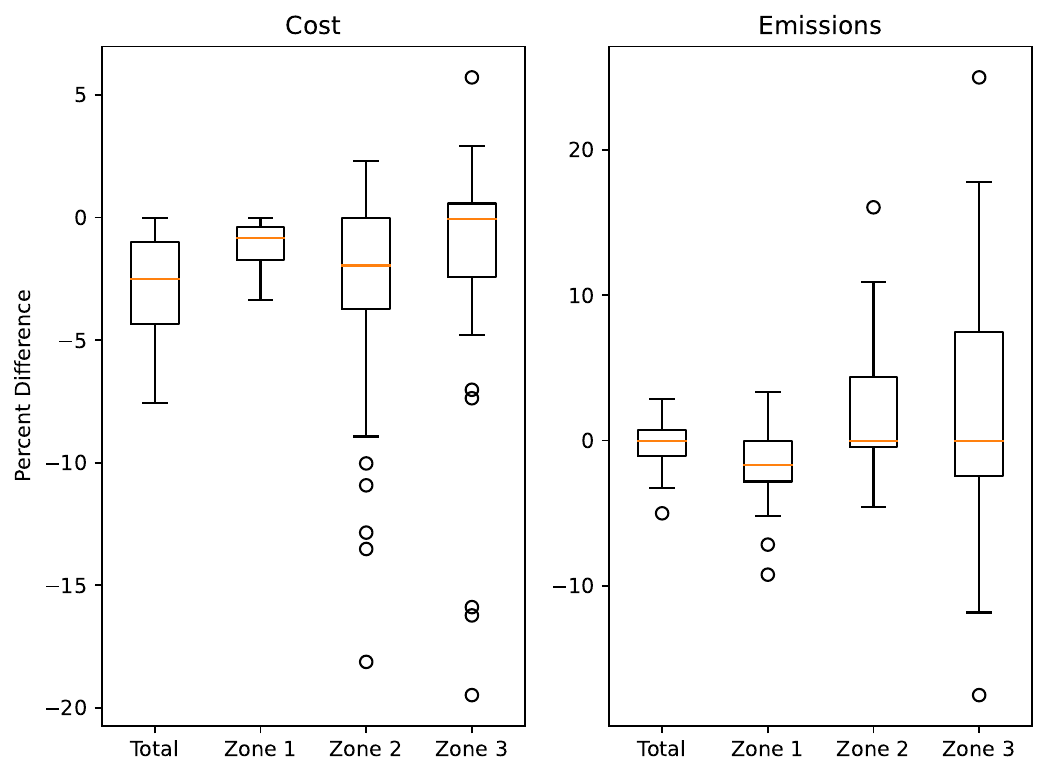}
    \caption{Boxplot of percent difference between metrics computed by MGCA and the full Capacity Expansion Model (n=50).}
    \label{fig:PercentDiff}
\end{figure}
We note several characteristics here. First, MGCA overestimates system cost by about 3\% on average and 7\% on the high end. This confirms that MGCA will always overestimate system cost due to the convex nature of the economic dispatch problem. It also confirms that the system cost is bounded by the budget constraint of the MGA problem (10\% for this problem). Note that because these are interior points, they are all below the original MGA budget constraint.  Second, MGCA has a harder time estimating zonal costs as they are not explicitly minimized in the CEM objective function (e.g. only total costs are minimized). We note here that Zone 3 has a much small overall system cost than Zones 1 and 2, making outliers appear much more extreme than for the other regions. Third, we note that emissions estimates are generally within 10\% error, but that several outliers exist. While, on average, these tend to underestimate emissions by 1-3\%, it is impossible to clearly bound these estimates beyond the overall bounds of the feasible space. Here again, zonal values are less reliable as they have greater variance relative to their magnitude. }
\paragraph{}{A second way to consider the accuracy of MGCA metric estimates is difference in zonal share as a proportion of the total of each metric. To calculate this we found each zone's percentage share of the total value of each metric for both the MGCA estimates and the full CEM solutions, then we took the difference of the MGCA and CEM values. The results of this are shown in Figure \ref{fig:PercentShare}. 
\begin{figure}[H]
    \centering
    \includegraphics[width=\linewidth]{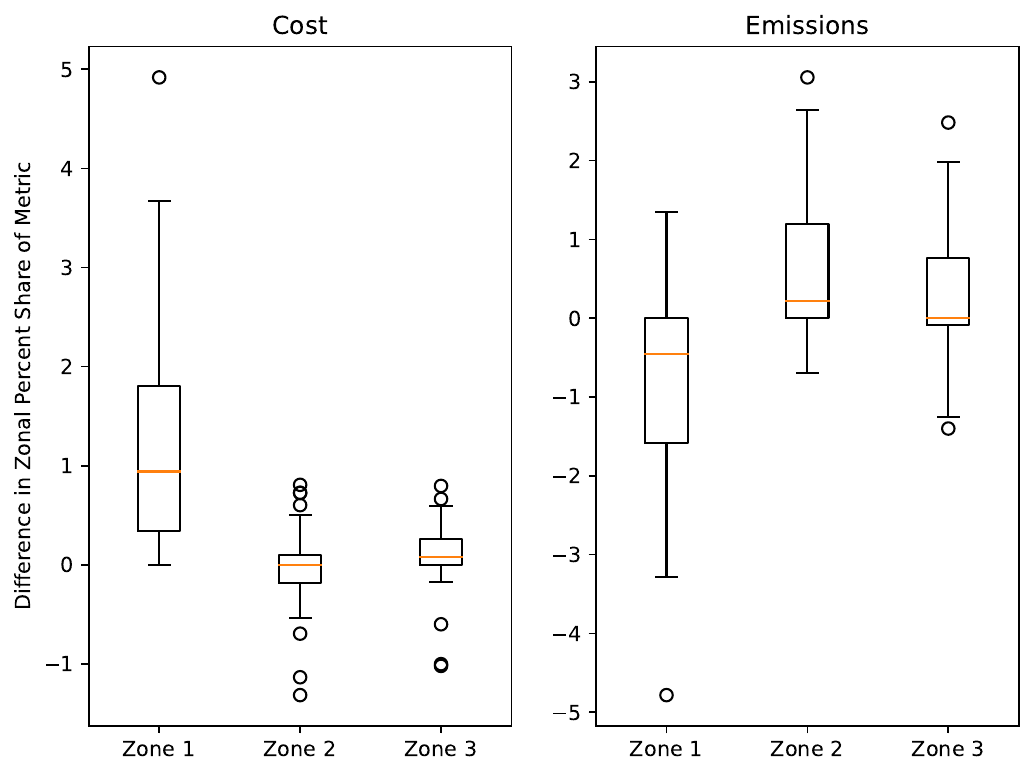}
    \caption{Boxplot of the difference in each zone's percentage share as calculated by MGCA and GenX with linearized unit commitment and economic dispatch (n = 50). }
    \label{fig:PercentShare}
\end{figure}
Based on these results, MGCA appears to get the share of the cost and emissions burden in each zone within an error range of 2-4\% in either direction. Thus the distribution of each metric is directionally correct. As discussed, zonal metrics have inherently greater variance as they are not explicitly optimized for, thus it is worth quantifying this level of uncertainty for each problem at hand. There is room for future work on further exploration of the reliability of operational metrics in MGCA interpolates and potential methods to provide firm bounds on differences between feasible solutions generated by MGCA and least-cost operational decisions.}
\section{Conclusions} \label{Conclusions}
\paragraph{}{In this paper, we introduced the novel MGA post-processing method Modelling to Generate Continuous Alternatives. Throughout Section \ref{Results}, we demonstrated several new capabilities unlocked through utilization of the properties of convex combinations, specifically the ability to evaluate MGA results under different budgets, the ability to incorporate additional constraints on the capacity and metric space, the ability to create approximate Pareto frontiers in post-processing, the ability to identify optimal capacity decisions for alternative objectives within the MGA space, and the ability to export and evaluate least-cost dispatch for points of interest.}
\paragraph{}{The key contributions of this work are fivefold: 
\begin{enumerate}
\item \textit{Dimensionality reduction to capacity decisions and metrics:} We introduced the concept of dimensionality reduction for MGA results to increase computation speed in post-processing while maintaining key information. We achieve this by retaining only capacity decisions, and metric values, which retains the strategic decisions of interest to MGA in the context of energy systems planning problems and their resulting outcomes of interest.
\item     \textit{Interpolating between MGA iterates through convex combinations}:  Pedersen et al. (2021) demonstrated MGA solutions could be viewed as a convex space with new points created through the convex combinations of the exterior points of simplicies within the space. We extend this result to use the MGA solution set as the set of exterior points for convex combinations, which provides greater flexibility in exploring the full range of the near-optimal feasible region of a problem (at least as such a space is approximated by previously computed MGA iterates) and avoids the computational obstacles of Delauney triangulation.  We used this capability to both generate capacities for new interpolates at random and to evaluate capacity mixes for lower budget slacks all computed in fractions of a second.
    \item \textit{Convex hull exploratory problem}: We incorporated MGA solutions as vertices in a convex combination-based optimization problem. Doing so allows significant additional capability to find and generate new interior solutions. We demonstrated a few uses of the feasibility problem by imposing constraints on the MGA capacity and metric space, conducting an additional MGA exploration with a tighter budget constraint, and approximating a Pareto frontier between cost and emissions.
    \item \textit{Convex combinations of metrics:} We showed that affine metrics of capacity variables could be evaluated for interpolated solutions through convex combinations. We also explored the behavior of operational metrics in reduced dimension interpolates, finding that the interpolates are guaranteed feasible and typically within 5-10\% of their actual value. 
    \item \textit{Exporting interpolated MGCA solutions for full analysis}: We demonstrate that MGCA solutions can be exported to a Capacity Expansion Model for full re-analysis and examine the accuracy of operational metric interpolations relative to economic dispatch in an 8760-hour capacity expansion model baseline. We find that projected operational metric values are always feasible in the original problem and are thus bounded by the extremes of the space. Total system cost, or whichever metric is the objective of the operational problem, will always be overestimated by MGCA, by a value that is less than or equal to the slack of the original problem's budget constraint. We provide an easy means for future users to export and resolve MGCA interpolates for economic dispatch in CEMs, allowing precise operational metric values to be calculated for any interpolate.
\end{enumerate}}
\subsection{Future Applications}
\paragraph{}{Overall, we demonstrate that MGCA is broadly applicable as a post-processing algorithm for Capacity Expansion Models using MGA. We have demonstrated that MGCA can allow users to explore MGA results in ways traditional MGA post-processing cannot, through generating new interior points with capacity metrics evaluated and operations metrics estimated, enabling spatial search for optimal values to affine functions of capacity, interpolating down to lower budgets, imposing new constraints or objectives and/or creating Pareto frontiers, all while providing a method to verify the operational characteristics of interpolates. Due to the dimensionality reduction methods introduced, all of this is carried out in small convex combination problems which solely require linear algebra, and in exploratory optimization problems at least 3-5 orders of magnitude smaller than the original CEMs. The resulting problems solve in fractions of a second, even using open-source solvers like GLPK, as used in this study.  This permits potential integration into a "live" interactive tool for exploration of the near-optimal feasible space where users (e.g. stakeholders or decisions makers) themselves can direct the exploration of alternative solutions, priorities and tradeoffs. Thus, further research should explore the incorporation of MGCA in live, interactive dashboards and decision support settings. Additionally, given the location-specific local opposition faced by energy infrastructure development around the world, research is needed on new down-scaling methods to enable live, user-malleable geospatial down-scaling of capacity siting decisions consistent with MGCA interpolates. All methods within this paper are extensible to linear programs with multiple solutions generally, though the discussion of dimensionality reduction and operational metric evaluation is specifically applicable to problems with well-defined strategic and operational decisions.}
\printbibliography
\end{document}